\documentclass[11pt]{article}

\usepackage{amsmath}
\usepackage{amsfonts}
\usepackage{amssymb}
\usepackage{amsthm}
\usepackage{enumitem}
\usepackage{mathtools}
\usepackage{setspace}
\usepackage{etoolbox}

\newtheoremstyle{mystyle}
{11pt}				
{11pt}				
{}					
{}					
{\bfseries}			
{}					
{5.5pt}				
{}					

\theoremstyle{mystyle}
\newtheorem{theorem}{Theorem}[section]

\newtheorem{proposition}[theorem]{Proposition}
\newtheorem{corollary}[theorem]{Corollary}

\setitemize{noitemsep, topsep=5.5pt, parsep=5.5pt, partopsep=0pt}

\allowdisplaybreaks
\appto\normalsize{
	\abovedisplayskip=5.5pt plus 2pt minus 2pt
	\belowdisplayskip=5.5pt plus 2pt minus 2pt
	\abovedisplayshortskip=5.5pt plus 2pt minus 2pt
	\belowdisplayshortskip=5.5pt plus 2pt minus 2pt}
\appto\small{
	\abovedisplayskip=5.5pt plus 2pt minus 2pt
	\belowdisplayskip=5.5pt plus 2pt minus 2pt
	\abovedisplayshortskip=5.5pt plus 2pt minus 2pt
	\belowdisplayshortskip=5.5pt plus 2pt minus 2pt}

\newcommand{\gap}{\vspace{11pt}}

\newcommand{\Aut}{\operatorname{Aut}}

\newcommand{\C}{\mathcal{C}}
\newcommand{\R}{\mathcal{R}}
\newcommand{\Rn}{\mathcal{R}^n}
\newcommand{\Rnp}{\mathcal{R}_+^n}
\newcommand{\Sn}{\mathcal{S}^n}

\newcommand{\Hn}{\mathcal{H}^n}
\newcommand{\V}{{\cal V}}

\newcommand{\spn}{\operatorname{span}}

\title{\bf On the connectedness of spectral sets and irreducibility of  spectral cones in Euclidean Jordan algebras}
\author{
 M. Seetharama Gowda\\
        Department of Mathematics and Statistics\\
        University of Maryland, Baltimore County\\
        Baltimore, Maryland 21250, USA\\
        gowda@umbc.edu\\[11pt]
        and\\[11pt]
	Juyoung Jeong\\
	Department of Mathematics and Statistics\\
	University of Maryland, Baltimore County\\
	Baltimore, Maryland 21250, USA\\
	juyoung1@umbc.edu
}

\date{\today}

\begin{document}

\maketitle

\begin{abstract}
Let $\V$ be a Euclidean Jordan algebra of rank $n$. A set $E$ in $\V$ is said to be a
{\it spectral set} if there exists a permutation invariant set
$Q$ in $\Rn$ such that $E=\lambda^{-1}(Q),$
where $\lambda:
\V\rightarrow \Rn$ is the {\it eigenvalue map} that takes  $x\in \V$ to $\lambda(x)$ (the vector of
eigenvalues of $x$ written in the decreasing order).
If the above $Q$ is also a convex cone, we say that $E$ is a {\it spectral cone}.
This paper deals with  connectedness and arcwise connectedness properties of spectral sets. By relying
on the result that in a simple Euclidean Jordan algebra, every eigenvalue orbit
$[x]:=\{y:\lambda(y)=\lambda(x)\}$
is arcwise connected, we show that
 if  a permutation invariant set $Q$ is connected (arcwise connected), then $\lambda^{-1}(Q)$
is connected (respectively, arcwise connected).
A related result is that in a simple Euclidean Jordan algebra, every pointed spectral cone is irreducible.

\end{abstract}

\vspace{1cm}
\\

\newpage


\section{Introduction}
Let $\V$ be a Euclidean Jordan algebra of rank $n$ \cite{faraut-koranyi}. The {\it eigenvalue map} $\lambda:
\V\rightarrow \Rn$ takes  $x\in \V$ to $\lambda(x)$, the vector of eigenvalues of
$x$ written in the decreasing order. 
A set $E$ in $\V$ is said to be a
{\it spectral set} if there exists a permutation invariant set
$Q$ in $\Rn$ such that $E=\lambda^{-1}(Q).$
If the above $Q$ is also a convex cone, we say that $E$ is a {\it spectral cone}.
A function $F:\V\rightarrow \R$ is said to be a {\it spectral function} if it is of the form 
$F=f\circ \lambda$, where $f:\Rn\rightarrow \R$ is a permutation invariant function. 
The above concepts are generalizations of similar concepts that have been studied in the 
settings of Euclidean $n$-space $\Rn$ and  
$\Sn$ ($\Hn$), the space of all $n \times n$ real (respectively, complex) Hermitian matrices, see for example,  
\cite{borwein-lewis}, \cite{chandrasekaran et al},  \cite{daniilidis et al},  \cite{daniilidis et al 2}, \cite{henrion-seeger}, \cite{iusem-seeger}, \cite{lewis}, \cite{lewis2}, \cite{lewis-sendov}, \cite{seeger}, and the references therein. In the case of $\Rn$, 
spectral sets/cones/functions are related to permutation invariance, and in $\Sn$ ($\Hn$) they  are precisely those that are invariant under linear transformations of the form
$X \rightarrow UXU^*$, where $U\in \R^{n\times n}$ is an  orthogonal (respectively, unitary) matrix.
In the general  setting of 
Euclidean Jordan algebras, they have been studied in several works, see 
\cite{baes}, \cite{gowda-jeong}, \cite{jeong-thesis}, \cite{jeong-gowda-spectral cone}, \cite{jeong-gowda-spectral set},
\cite{ramirez et al}, \cite{sossa}, and \cite{sun-sun}.

Focusing on topological/convexity/linearity properties of spectral sets/cones, in two recent papers 
Jeong and Gowda \cite{jeong-gowda-spectral cone}, \cite{jeong-gowda-spectral set}   
show that the multivalued map $\lambda^{-1}$ from $\Rn$ to $\V$ behaves like a  
linear isomorphism on certain types of permutation invariant sets. Specifically, the following results are shown (where part of the first result is due to Baes \cite{baes}):

\begin{proposition} \label{basic properties}
{\it Let $Q$, $Q_1$, and $Q_2$ be permutation invariant sets in $\Rn$ with $Q_1$ and $Q_2$ convex. Let $\alpha\in \R$.
Then
\begin{itemize}
\item [(a)] $\lambda^{-1}(Q)$ is open/closed/compact/convex/cone in $\V$ if and only if $Q$ is so in $\Rn$.  
Moreover,
$\lambda^{-1}(Q^\#)=[\lambda^{-1}(Q)]^\#$, where $\#$ 
denotes any operation of taking closure, interior, boundary, or convex/conic hull. 
\item [(b)]
$\lambda^{-1}(Q_1+Q_2)=\lambda^{-1}(Q_1)+\lambda^{-1}(Q_2)\,\,\mbox{and}\,\, \lambda^{-1}(\alpha Q_1)=\alpha \lambda^{-1}(Q_1).$
\item [(c)] When $Q$ is a convex cone,  $\lambda^{-1}(Q)$ is pointed if and only if $Q$  is pointed.
\end{itemize}
}
\end{proposition}

Motivated by these results, 
we ask if connectedness and  arcwise (=pathwise) connected properties of permutation invariant $Q$
carry over to $\lambda^{-1}(Q)$. In this paper, we answer these affirmatively by relying on the  result that in a 
simple Euclidean Jordan algebra, for any element $x$, the eigenvalue orbit 
$[x]:=\{y: \lambda(y)=\lambda(x)\}$
is  arcwise connected. This result will also be used to show that in any simple Euclidean Jordan algebra, every
pointed spectral cone is irreducible.  
\section{Preliminaries}
We let $\Rn$ denote the Euclidean $n$-space (where the vectors are regarded as column vectors or row vectors depending on the context). 
In $\Rn$, we denote the standard coordinate vectors by $c_1,c_2,\ldots, c_n$, where 
$c_i$ is the vector with one in the $i$th slot and zeros elsewhere. We use the notation $\Sigma_n$ to denote the set of all $n\times n$ permutation matrices. For any set $S$ in $\Rn$, we let $\Sigma_n(S):=\{\sigma(s):\sigma\in \Sigma_n, s\in S\}.$
For any vector $q\in \Rn$ and a set $Q\subseteq \Rn$, $q^\downarrow$ denotes the decreasing rearrangement of $x$ 
(i.e., the entries of $q^\downarrow$ satisfy $q^\downarrow_1\geq q^\downarrow_2\geq \cdots\geq q^\downarrow_n$) and 
$Q^{\downarrow}:=\{q^\downarrow:q\in Q\}$. 
A non-empty set $Q$ in $\Rn$ is said to be {\it permutation invariant} if $\sigma(Q)=Q$ for all
$\sigma\in \Sigma_n$. \\
We assume that the reader is familiar with standard topological notions/results dealing with  (arcwise) connectedness, components, etc. Recall that a set 
$S$ (say, in a topological space) is {\it connected} if it is not the union of two nonempty separated sets (where separation takes the form $\overline{A}\cap B=\emptyset=A\cap \overline{B}$, with  the `overline' indicating the closure) \cite{rudin}. The set $S$ is arcwise (= pathwise) connected if  any two points of $S$ can be joined by a continuous arc ($=$ continuous image of an interval in $\R$) that lies inside $S$. Connected (arcwise connected) components of a set $S$ are maximal connected (respectively, arcwise connected) sets in $S$.
\\
For  basic things related to Euclidean Jordan algebras, we refer to \cite{faraut-koranyi} or \cite{jeong-gowda-spectral set} for a summary. 
{\it Throughout this paper, $\V$ denotes a Euclidean Jordan algebra of rank $n$.} 
Recall that a  Euclidean Jordan algebra $\V$ is  {\it simple} if it is not a direct product of  nonzero Euclidean Jordan
algebras (or equivalently, if it does not contain any non-trivial ideal). It is known (see \cite{faraut-koranyi}, Prop. III.4.4)
that any nonzero Euclidean Jordan algebra is, in a unique way, a direct product/sum of simple Euclidean Jordan algebras.
Moreover, there are (in the isomorphic sense) only five simple algebras, two of which are:
$\mathcal{S}^n$, the algebra  of $n \times n$ real symmetric matrices, and $\mathcal{H}^n$, the algebra of $n \times n$ complex Hermitian matrices. The other three are: $n\times n$ quaternion Hermitian matrices, $3\times 3$ octonion Hermitian matrices, and the Jordan spin algebra.
In $\V$, each element $x$ has a  spectral decomposition: $x=q_1e_1+q_2e_2+\ldots+q_ne_n$, 
where $\{e_1,e_2,\ldots,e_n\}$ is a Jordan frame in $\V$ and the real numbers 
$q_1,q_2,\ldots, q_n$ are (called) the eigenvalues of $x$. 
Then $\lambda(x)=\big (\lambda_1(x),\lambda_2(x),\ldots,\lambda_n(x)\big )$ denotes the vector of eigenvalues of $x$ written in the decreasing order. We note that
\begin{center}
$\lambda:\V\rightarrow \Rn$ is continuous (\cite{baes}, Corollary 24) and $\lambda(\V)=(\Rn)^\downarrow.$
\end{center}
For any $x\in \V$, we let
$$[x]:=\{y\in \V: \lambda(y)=\lambda(x)\}$$
denote the `eigenvalue orbit' of $x$. (The notation $[x]_{\V}$ will be used when  more than one algebra is involved.) 
Suppose $\V$ is a Cartesian product (or a direct sum)
$\V=\V^{(1)}\times \V^{(2)}\times \cdots\times \V^{(N)}, $ where each $\V^{(i)}$ is simple. 
It is easy to see that in  $\V$,   
an element $c$ is a primitive idempotent if and only if it is of the form $(0,0,\ldots,0, c_i,0,\ldots,0)$, where $c_i$ is a primitive idempotent in $\V^{(i)}$ for some $i$. Applying this to elements of a Jordan frame in $\V$, we see that 
the eigenvalues of any $x=(x^{(1)},x^{(2)},\ldots, x^{(N)})\in \V$ comprise of the eigenvalues of $x^{(i)}$ in $\V^{(i)}$, $i=1,2,\ldots, N$. For such an $x$,  we define the `restricted eigenvalue orbit' by
$$[x]_r:=\Big \{ y=\Big(y^{(1)},y^{(2)},\ldots, y^{(N)}\Big): y^{(i)}\in [x^{(i)}]_{\V^{(i)}}\,\,\mbox{for all}\,\,i\Big \}.$$
 We note that $[x]_r\subseteq [x]$ with equality when $\V$ is simple.
To see an example, let $\V=\Rn=\R\times \R\times \cdots\times\R$ ($n>1)$. 
Recalling that $\{c_1,c_2,\ldots, c_n\}$ denotes the set of standard
coordinate vectors in $\Rn$, we see that $[c_1]=\{c_1,c_2,\ldots, c_n\}$, while $[c_1]_r=\{c_1\}$. 
\\
Given a Jordan frame ${\cal E}=\{e_1,e_2,\ldots,e_n\}$ in $\V$, we assume that its listing/enumeration is fixed and 
 define, for any $q=(q_1,q_2,\ldots, q_n)\in \Rn$,
$$q*{\cal E}:=q_1e_1+q_2e_2+\cdots+q_ne_n.$$ 
We note that 
$$\lambda(q*{\cal E})=q^\downarrow$$
and when ${\cal E}$ is fixed,  
$\Theta:q\mapsto q*{\cal E}$ is a continuous map from $\Rn$ to $\V$. 

 Recall that a nonempty set in $\V$  is a {\it spectral set} if it is of the form $\lambda^{-1}(Q)$ for some permutation invariant set $Q$ in $\Rn$. 
We will freely use the results in the following (easily verifiable)  proposition.

\begin{proposition}\label{elementary properties of lambda inverse}
 {\it Let $P$ and $Q$ be nonempty subsets of $\Rn$ with $Q$  permutation invariant. Then, 
\begin{itemize}
\item [$\bullet$] $Q=\Sigma_n(Q^\downarrow).$
\item [$\bullet$] $x\in \lambda^{-1}(Q)$ if and only if $x=q*{\cal E}$ for some Jordan frame ${\cal E}$ and $q\in Q$.
\item [$\bullet$] $\lambda\big ( \lambda^{-1}(P)\big )= P\cap P^\downarrow$ and $\lambda^{-1}(P)=\lambda^{-1}(P\cap P^\downarrow)$. In particular, for any $\Omega\subseteq (\Rn)^\downarrow$, $\lambda(\lambda^{-1}(\Omega))=\Omega$ and $\lambda^{-1}(Q)=\lambda^{-1}(Q^\downarrow).$ 
\item [$\bullet$] Spectral sets can be generated by taking $\lambda-$inverse images of subsets of $(\Rn)^\downarrow$. In fact, any set of the form $\lambda^{-1}(P)$ is a spectral set. 
\item [$\bullet$] 
The correspondence $\Omega\mapsto \lambda^{-1}(\Omega)$ is  one-to-one and onto  between nonempty subsets of $(\Rn)^\downarrow$ and spectral sets in $\V$.
\item [$\bullet$] For a Jordan frame ${\cal E}=\{e_1,e_2,\ldots,e_n\}$, the set
$\Theta(Q)=\{q*{\cal E}:q\in Q\}$ is independent of the listing of elements in 
${\cal E}$.
\item [$\bullet$] $\lambda^{-1}(Q)$ is the union of sets of the form 
$\Theta(Q)$ as ${\cal E}$ varies over all Jordan frames.
\end{itemize}
}\end{proposition}
Recall that an {\it automorphism} $\phi$ of $\V$ is an invertible linear transformation on $\V$ that satisfies the condition $\phi(x\circ y)=\phi(x)\circ \phi(y)$ for all $x,y\in \V$. {\it The set of all such transformations is denoted by $\Aut(\V)$ and we let
 $G$ denote the connected component of the identity transformation (henceforth called the  component of identity) 
in $\Aut(\V)$.}
For example, $\Aut(\mathcal{S}^n$) consists of transformations of the form
$\phi(X):=UXU^T\,\,(X\in \mathcal{S}^n) $ where $U\in \R^{n\times n}$ is an orthogonal matrix. When $U$ has determinant one, 
such a $\phi$ belongs to the corresponding $G$. Also, $\Aut(\mathcal{H}^n$) consists of transformations of the form
$\phi(X):=UXU^*\,\,(X\in \mathcal{H}^n),$ where $U$ is a unitary matrix. In this case, $G=\Aut(\mathcal{H}^n)$ (\cite{hall}, page 15).
\\
We need the following result connecting eigenvalue orbits and spectral sets.

\begin{proposition}\label{spectral set properties}
{\it \begin{itemize}
\item [(i)] The following inclusions hold:
$$\big \{\phi(x):\phi\in G\big\}\subseteq \big \{\phi(x):\phi\in \Aut(\V)\big \}\subseteq \big [x\big ]\quad(x\in \V).$$
These become equalities when $\V$ is simple.
\item [(ii)] A set $E$ in $\V$ is a spectral set if and only if for all $x\in E$, $[x]\subseteq E$.
\item [(iii)] If $E$ is a spectral set in $\V$, then for all $\phi\in \Aut(\V)$, $\phi(E)=E$. Converse holds when $\V$ is simple or $\Rn$.
\end{itemize}
}
\end{proposition}

\noindent{\bf Proof.} $(i)$ As automorphisms preserve Jordan frames and eigenvalues, the stated inclusions follow. Now suppose
$\V$ is simple and let $y\in [x]$ so that $\lambda(x)=\lambda(y)$.
We write the spectral decompositions $x=\sum_{1}^{n}\lambda_i(x)e_i$ and $y=\sum_{1}^{n}\lambda_i(y)f_i$,
where $\{e_1,e_2,\ldots, e_n\}$ and $\{f_1,f_2,\ldots, f_n\}$ are Jordan frames in $\V$. Since $\V$ is simple, by
Corollary IV.2.7 in \cite{faraut-koranyi}, there is an automorphism in $G$ that takes one Jordan frame to the other.
Hence, we may write $y=\phi(x)$ for some $\phi\in G$. Thus, by $(i)$, $[x]=\{\phi(x):\phi\in G\}$.\\
$(ii)$ This is given in Theorem 1, \cite{jeong-gowda-spectral set}.
\\
$(iii)$ This is given in Theorem 2, \cite{jeong-gowda-spectral set}.
$\hfill$ $\qed$
\\

In view of Item $(iii)$ in the above theorem, we see that in $\mathcal{S}^n$, a set $E$ is spectral if and only if 
$$X\in E\Rightarrow UXU^T\in E\,\,\mbox{for all orthogonal matrices}\,\,U\in \R^{n\times n}.$$

\section{Main connectedness results}
To motivate our discussion and results, we start with two examples.
\\

\noindent{\bf Example (1)} Let $n>1$ and consider  the non-simple algebra $\V=\Rn$. Then, $\Aut(\Rn)=\Sigma_n$ and $G$ consists of just the identity matrix. 
For any $q\in \Rn$, $\lambda(q)=q^\downarrow$.  
Also, for any permutation invariant set $Q$ in $\Rn$, $\lambda^{-1}(Q)=Q$. 
Thus, in this setting, if a permutation invariant set is connected (arcwise connected), then 
its $\lambda$-inverse image is connected (respectively, arcwise connected). Can we improve this by assuming only 
the connectedness (arcwise connectedness) of $Q^\downarrow$? To answer this, let $Q:=\{c_1,c_2,\ldots, c_n\}$
denote the set of standard coordinate vectors in $\Rn$. Then $Q$ is permutation invariant and $Q^\downarrow =\{c_1\}$. We see that 
while $Q^\downarrow$ is connected (arcwise connected), $\lambda^{-1}(Q)$ is not connected.
\\

\noindent{\bf Example (2)} Let $n>1$ and consider the simple algebra $\V={\cal S}^n$.  
As noted previously,  $\Aut(\Sn$) consists of transformations of the form
$\phi(X):=UXU^T\,\,(X\in \mathcal{S}^n) $ where $U\in \R^{n\times n}$ is an orthogonal matrix. 
Now consider the set $Q:=\{c_1,c_2,\ldots, c_n\}$ of Example 1. As noted earlier, $Q$ is a permutation invariant set with 
$Q^\downarrow$ connected.  
With $D$ denoting the diagonal matrix with $c_1$ on its diagonal, we have
$$\lambda^{-1}(Q)=\lambda^{-1}(Q^\downarrow)=\Big \{UDU^T:\,\,U\in \R^{n\times n}\,\,\mbox{is orthogonal}\Big \}=\Big \{uu^T: \,u\in \Rn, ||u||=1\Big \}.$$ (The vector $u$ that appears above is the first column of $U$.) 
As $n>1$, the unit sphere in $\Rn$ is arcwise connected. Hence, $\lambda^{-1}(Q)$ is arcwise connected.
Thus, in contrast to Example 1, a weaker hypothesis involving connectedness (arcwise connectedness) of $Q^\downarrow$ 
suffices  
to get the connectedness (arcwise connectedness) of $\lambda^{-1}(Q)$.  
We show below that in the setting of a simple algebra, such a statement holds for any permutation invariant set. \\

The following is a basic  (possibly known) result   about the connectedness of eigenvalue orbits. Recall that $G$ is the connected component of identity in $\Aut(\V)$.

\begin{theorem}\label{G is arcwise connected}
{\it The following statements hold:
\begin{itemize}
\item [(a)] $G$ is arcwise connected. 
\item [(b)] When $\V$ is simple, for any $x\in \V$, $[x]$
is arcwise connected.
\item [(c)] For any $x\in \V$, $[x]_r$
is arcwise connected.
\end{itemize}
}
\end{theorem}

\noindent{\bf Proof.} $(a)$ As $\Aut(\V)$ is a (matrix) Lie group (see \cite{faraut-koranyi}, page 36), 
its connected component of identity, namely $G$,  is
arcwise connected, see \cite{duistermaat-kolk}, Theorem 1.9.1. \\
$(b)$ Now suppose $\V$ is simple and consider any $x\in \V$. Then, by Item $(i)$
in the previous proposition, $[x]=\{\phi(x):\phi\in G\}$. As $G$ is arcwise connected and the map $\phi\mapsto \phi(x)$ is continuous, $[x]$ is also arcwise connected.
\\
$(c)$ If $\V$ is simple, $[x]_r=[x]$ for any $x$. Then, the result follows from $(b)$. Suppose
$\V$ is non-simple, and let $\V$ be a product of simple algebras:
$\V=\V^{(1)}\times \V^{(2)}\times \cdots\times \V^{(N)}, $ where each 
$\V^{(i)}$ is simple. Let  $x=(x^{(1)},x^{(2)},\ldots, x^{(N)})\in \V$. 
Then, by $(b)$, each $[x^{(i)}]_{\V^{(i)}}$ is arcwise connected. Hence, the set
$$
[x]_r=[x^{(1)}]_{\V^{(1)}}\times [x^{(2)}]_{\V^{(2)}}\times \cdots\times [x^{(N)}]_{\V^{(N)}},$$
being a product of arcwise connected sets, is arcwise connected.
$\hfill$ $\qed$

\gap

To illustrate  Items $(b)$ and $(c)$ in the  above result, we let
 $x$  be a primitive idempotent in a simple algebra $\V$. Then,  $[x]$, which is
the set of all primitive idempotents, is arcwise connected. (This result is known, see
\cite{faraut-koranyi}, page 71.)  This may fail  if $\V$ is not simple:
 let $\V = \Rn$ ($n>1$) and $x=c_1$.
Then, $[x]=\{c_1,c_2,\dots, c_n\}$ is not connected.
However, $[x]_r=\{c_1\}$ (a singleton set) is arcwise connected.
\\

We now state a necessary condition for a $\lambda$-inverse image to be connected (arcwise connected).

\begin{proposition}\label{necessary condition}
{\it 
Let $P\subseteq \Rn$. If  $\lambda^{-1}(P)$ is connected (arcwise connected) in $\V$, then $P\cap P^\downarrow$ is connected (arcwise connected) in $\Rn$.
In particular, if $Q$ is permutation invariant and $\lambda^{-1}(Q)$ is connected (arcwise connected) in $\V$, then $Q^\downarrow$ is connected (arcwise connected) in $\Rn$.
}\end{proposition}

\noindent{\bf Proof.} The first statement follows from the continuity of the eigenvalue map $\lambda$ and 
the equality $P\cap P^\downarrow =\lambda(\lambda^{-1}(P))$. The second one follows from the equality 
$Q\cap Q^\downarrow=Q^\downarrow$.
$\hfill$ $\qed$
\\

The next two results deal with sufficient conditions.

\begin{theorem}\label{connectedness result}
{\it
Suppose  $Q\subseteq \Rn$ is  permutation invariant and one of the following conditions holds.
\begin{itemize}
\item [(a)] $\V$ is simple and $Q^\downarrow$ is connected (arcwise connected).
\item [(b)] $Q$ is connected (arcwise connected).
\end{itemize}
Then
$\lambda^{-1}(Q)$ is  connected (respectively, arcwise connected) in $\V$.
}
\end{theorem}

\noindent{\bf Proof.}
Suppose condition $(a)$ holds so that $\V$ is simple and  $Q^\downarrow$ is
 connected (arcwise connected).
We fix a Jordan frame ${\cal E}=\{e_1,e_2,\ldots, e_n\}$ and consider the
continuous map $\Theta:\Rn\rightarrow V$ defined by $\Theta(q):=q*{\cal E}$, see Section 2. Then, the image
$\Theta(Q^\downarrow)$, which is a subset of $\lambda^{-1}(Q^\downarrow)$, is connected (respectively, arcwise connected).
Now let $x\in \lambda^{-1}(Q^\downarrow)$ be arbitrary. Then, $\lambda(x)\in Q^\downarrow$ and there is a Jordan frame $\{f_1,f_2,\ldots, f_n\}$ such that
$x=\lambda_1(x)f_1+\lambda_2(x)f_2+\cdots+\lambda_n(x)f_n$.
As $\V$ is simple, from Theorem \ref{G is arcwise connected}, $[x]$ is arcwise connected. Also,
$$\lambda_1(x)e_1+\lambda_2(x)e_2+\cdots+\lambda_n(x)e_n\in \Theta(Q^\downarrow)\cap [x].$$
So the sets $\Theta(Q^\downarrow)$ and $[x]$ are connected (respectively, arcwise connected) and their intersection is nonempty.
It follows that their union is also connected (respectively, arcwise connected). As $x$ is arbitrary in $\lambda^{-1}(Q^\downarrow)$, the connected component (respectively, arcwise connected component) of $\lambda^{-1}(Q^\downarrow)$ that contains
$\Theta(Q^\downarrow)$ must be $\lambda^{-1}(Q^\downarrow)$ itself. This proves that $\lambda^{-1}(Q^\downarrow)$ is
connected (respectively, arcwise connected) under the condition $(a)$. The stated assertion follows since
$\lambda^{-1}(Q)=\lambda^{-1}(Q^\downarrow).$\\
Now suppose condition $(b)$ holds. If $\V$ is simple, we can use the previous argument as $Q^\downarrow $ (which is the image of $Q$ under the continuous 
map $q\mapsto q^\downarrow$) is connected (arcwise connected). So, assume that $\V$ is non-simple and write it as an orthogonal direct sum 
$\V=\V^{(1)}\oplus \V^{(2)}\oplus\cdots\oplus \V^{(N)}, $ where each $\V^{(i)}$ is simple.
In $\V$, we fix a   
 Jordan frame ${\cal E}:=\{e_1,e_2,\ldots, e_n\}$ and consider the continuous map $\Theta:\Rn\rightarrow \V$ defined by $\Theta(q)=q*{\cal E}$. 
As $Q$ is connected (arcwise connected) and $\Theta$  is continuous, $\Theta(Q)$ is connected (arcwise connected). Also, as $Q$ is permutation invariant, any rearrangement/listing of elements of ${\cal E}$ will not alter the set $\Theta(Q)$. Hence, we may assume without loss of generality that
${\cal E}=\bigcup_{1}^{N}{\cal E}^{(i)}$, where each ${\cal E}^{(i)}$ is a Jordan frame in $V^{(i)}$.
Writing 
$q\in \Rn$ in the block form $q=(q^{(1)},q^{(2)},\ldots, q^{(N)})$, we see that
 $q*{\cal E}=q^{(1)}*{\cal E}^{(1)}+q^{(2)}*{\cal E}^{(2)}+\cdots+q^{(N)}*{\cal E}^{(N)}.$ 
Now, let $x\in \lambda^{-1}(Q)$. Then, there exist a Jordan frame 
${\cal F}$ in $\V$ and a $p\in Q$  such that $x=p*{\cal F}$. Similar to the 
above, we may write  $x=p*{\cal F}=p^{(1)}*{\cal F}^{(1)}+p^{(2)}*{\cal F}^{(2)}+\cdots+p^{(N)}*{\cal F}^{(N)},$
where ${\cal F}^{(i)}$ is a Jordan frame in $V^{(i)}$.
Let $y:= p*{\cal E}=p^{(1)}*{\cal E}^{(1)}+p^{(2)}*{\cal E}^{(2)}+\cdots+p^{(N)}*{\cal E}^{(N)}.$ As 
$\lambda\big(p^{(i)}*{\cal E}^{(i)}\big)=\lambda\big(p^{(i)}*{\cal F}^{(i)}\big)$ for all $i$, it follows that $y\in [x]_r$.
Hence, $y\in \Theta(Q)\cap [x]_r$. As $\Theta(Q)$ is connected (arcwise connected) and $[x]_r$ is arcwise connected (by Theorem \ref{G is arcwise connected}), it follows that $\Theta(Q)\cup [x]_r$ is connected (respectively, arcwise connected). 
Since $x$ is arbitrary in $\lambda^{-1}(Q)$, the connected (arcwise connected) component of $\lambda^{-1}(Q)$ must be $\lambda^{-1}(Q)$. This proves that 
$\lambda^{-1}(Q)$ is connected (respectively, arcwise connected).
$\hfill$ $\qed$
 
\gap

\begin{corollary} \label{corollary}
{\it When $\V$ is simple, the following statements hold:
\begin{itemize}
\item [(i)] If $\Omega$ is any subset of $(\Rn)^\downarrow$ that is connected (arcwise connected), then $\lambda^{-1}(\Omega)$ is connected (respectively, arcwise connected).
\item [(ii)] If $P$ is any subset of $\Rn$ such that $P\cap P^\downarrow$ is connected (arcwise connected), then 
$\lambda^{-1}(P)$ is connected (respectively, arcwise connected).
\end{itemize}
}\end{corollary}

\noindent{\bf Proof.} $(i)$ Let $\Omega$  be connected (arcwise connected) subset of $(\Rn)^\downarrow$. 
Then $Q:=\Sigma_n(\Omega)$ is permutation invariant and $Q^\downarrow =\Omega$. Hence condition $(a)$ in the above 
theorem applies. Since  $\lambda^{-1}(\Omega)=\lambda^{-1}(Q)$, we see that 
$\lambda^{-1}(\Omega)$  is connected (arcwise connected).\\
$(ii)$ Suppose $P$ is any subset of $\Rn$ such that $P\cap P^\downarrow$ is connected (arcwise connected). As 
$P\cap P^\downarrow$ is a subset of   $(\Rn)^\downarrow$ and $\lambda^{-1}(P)=\lambda^{-1}(P\cap P^\downarrow),$ the stated result follows from $(i)$ applied to $\Omega:=P\cap P^\downarrow$.
$\hfill$ $\qed$

\gap

\noindent{\bf Remarks (1)} Suppose $\V$ is simple. In view of Proposition \ref{necessary condition}, for any set $P$ in $\Rn$, connectedness (arcwise connectedness) of $P\cap P^\downarrow$ is  necessary and sufficient for that of $\lambda^{-1}(P)$. Thus, {\it a spectral set $E$ in a simple algebra $\V$ is connected 
(arcwise connected) if and only if $\lambda(E)$ is connected (arcwise connected) in $\Rn$.} \\

\noindent {\bf (2)} The conclusions in the above corollary may fail if $\V$ is not simple. For example in $\R^3$, $c_1=(1,0,0)$ is on the boundary of 
$(\R^3)^\downarrow$. Letting $\V=\R^3$ and $\Omega=\{c_1\}$, we see that
$\lambda^{-1}(\Omega)=\{c_1,c_2,c_3\}$ (set of standard coordinate vectors in $\R^3$) is not connected. In the same setting, it is easy to 
 construct an arcwise connected set
$P$ such that $P\cap P^\downarrow =\{c_1\}$. (For example, $P$ could be a 
circle through $c_1$ such that $P\backslash \{c_1\}$ is outside of $(\R^3)^\downarrow$.) For such a $P$, $\lambda^{-1}(P)=\lambda^{-1}(P\cap P^\downarrow)=\{c_1,c_2,c_3\}$ is not connected. 
\\

\noindent{\bf (3)} Motivated by the above two results, one may ask if $\lambda$-inverse image of a simply connected permutation invariant set $Q$
in $\Rn$ is simply connected in $\V$. While the answer to this is not clear, we note  
that in ${\cal S}^2$, the set of all
primitive idempotents is given by (\cite{faraut-koranyi}, page 71)
$$ \Big \{\left [ \begin{array}{cc}
\cos^2\theta & \cos\theta\sin\theta\\
\cos\theta\sin\theta & \sin^2\theta\end{array}\right ],\,0\leq \theta\leq \pi\Big \}.$$
This is homeomorphic to a circle, hence not simply connected. However, it is of the form $\lambda^{-1}(Q^\downarrow)$, where
$Q$ is the set of standard coordinate vectors in $\R^2$. While this $Q$ is not simply connected, $Q^\downarrow$,
 consisting of only one element, is simply connected.
So the counterpart of $(a)$ in Theorem \ref{connectedness result} may not hold for simple connectedness.

\gap

We end this section by mentioning a classical result of Ky Fan. Suppose $Q$ is a permutation invariant set that 
satisfies one of the conditions of Theorem \ref{connectedness result}. Assume further that $Q$ is compact. 
Then, $\lambda^{-1}(Q)$ is connected and 
(by Proposition \ref{basic properties}) compact in $\V$. Hence, for any $c\in \V$, the image of $\lambda^{-1}(Q)$ under the continuous function $x\mapsto \langle c,x\rangle$ is connected and compact in $\R$. So, there exists real numbers $\delta$ and $\Delta$ such that 
\begin{equation}\label{fan type equality}
\Big \{\langle c,x\rangle:x\in \lambda^{-1}(Q)\Big \}=[\delta,\Delta].
\end{equation}
(With some additional work, it is possible to describe the forms of $\delta$ and $\Delta$.)
To see a special case, consider two  matrices $C$ and $A$ in $\Hn$ and let 
$Q=\Sigma_n(\{\lambda(A)\})$. Clearly, $Q$ is permutation invariant and compact (actually, finite). Since 
 $\Hn$ is simple and $Q^\downarrow$ (a singleton set) is connected, by Theorem \ref{connectedness result} (or by Theorem \ref{G is arcwise connected}), $\lambda^{-1}(Q)=[A]=\{UAU^*:\,\,U\in \C^{n\times n}\,\,\mbox{is unitary}\}$ is connected in $\Hn$.  As $\langle X,Y\rangle=tr(XY)$ in $\Hn$, (\ref{fan type equality}) reads:
$$\Big \{tr(CUAU^*):U\in \C^{n\times n}\,\,\mbox{is unitary}\Big \}=[\delta,\Delta].$$ 
With  $\Delta:=\langle \lambda(C),\lambda(A)\rangle$ and $\delta:=\langle \lambda^{\uparrow}(C),\lambda(A)\rangle$, where 
$\lambda^{\uparrow}(C)$ denotes the increasing rearrangement of $\lambda(C)$, this statement is due to  Fan \cite{fan} 
(see also \cite{tam}, Corollary 1.6).  

\section{Components of a spectral set}

We now describe connected (arcwise connected) components of a spectral set in a simple algebra.

\begin{theorem}\label{components of a spectral set}
{\it Let $\V$ be  simple and $E:=\lambda^{-1}(Q)$ be a spectral set in $\V$, where $Q$ is permutation invariant in $\Rn$. 
Then, every connected (arcwise connected) component of $E$ is a spectral set. 
 Moreover, $C\rightarrow \lambda(C)$  is a  one-to-one correspondence between connected (arcwise connected) components of $E$ and those of $Q^\downarrow$.
}
\end{theorem}

\noindent{\bf Proof.}
We consider the case of connected components. The case of arcwise connected components is similar.
Let $C$ be a connected component of $E$. As $\V$ is simple, for any $x\in C$, $[x]$ is (arcwise) connected (by Theorem 
\ref{G is arcwise connected})
and $[x]\cap C$ contains $x$. Hence, $[x]\subseteq C$. By Proposition \ref{spectral set properties}(ii),
$C$ is a spectral set; we may now write $C$ as $C=\lambda^{-1}(P)$, where $P$ is a permutation invariant set in $\Rn$. As $C$ is connected, its image 
$\lambda(C)=P^\downarrow$  is connected. We claim that this set is a connected component of $Q^\downarrow.$ 
To simplify the notation, let $\Omega:=P^\downarrow.$ Then,
$$C\subseteq E\Rightarrow P^\downarrow =\lambda(C)\subseteq \lambda(E)= Q^\downarrow\Rightarrow
\Omega \subseteq Q^\downarrow.$$
So, $\Omega$ is a connected subset of $Q^\downarrow$. 
Let $\Omega^*$ be the connected component of $Q^\downarrow$ that contains $\Omega$, so 
$$\Omega\subseteq \Omega^*\subseteq Q^\downarrow.$$
Then, 
$$C=\lambda^{-1}(\Omega)\subseteq \lambda^{-1}(\Omega^*)\subseteq \lambda^{-1}(Q^\downarrow)=\lambda^{-1}(Q)=E.$$
Since $\V$ is simple, by Item $(i)$ in Corollary \ref{corollary}, $\lambda^{-1}(\Omega^*)$ is connected.  By our 
assumption that $C$ is a connected component of $E$ we get $C=\lambda^{-1}(\Omega^*)$, leading to 
$\Omega=\lambda(C)=\lambda((\lambda^{-1}(\Omega^*))=\Omega^*,$ the last equality is due to Proposition \ref{elementary properties of lambda inverse}. 
Hence, $\Omega$  is a connected component of $Q^\downarrow.$
\\Now we show that every connected component of $Q^\downarrow$ arises this way. Let $\Omega$ be a connected  component
of $Q^\downarrow$. By Item $(i)$ in Corollary \ref{corollary}, $C:=\lambda^{-1}(\Omega)$ is connected in $E$. Suppose $C^*$ is the connected component of $E$ that contains $C$. We show that $C=C^*$. Now, 
$$C\subseteq C^*\subseteq E\Rightarrow \lambda(C)\subseteq \lambda(C^*)\subseteq \lambda(E).$$
We have $\lambda(C)=\lambda(\lambda^{-1}(\Omega))=\Omega$, where the second equality is due to  Proposition 
\ref{elementary properties of lambda inverse}.  It follows that 
$\Omega\subseteq \lambda(C^*)\subseteq Q^\downarrow$. As $\Omega$ is a connected component and $\lambda(C^*)$ is connected,
we must have $\Omega=\lambda(C^*)$, that is, $\Omega=\lambda(C)=\lambda(C^*)$. 
So, $C=\lambda^{-1}(\Omega)=\lambda^{-1}(\lambda(C^*))\supseteq C^*$. Since $C\subseteq C^*$, we conclude that $C=C^*$ and that
$C$ is a connected component of $E$.  
Finally, the concluding statement in the theorem about one-to-one and onto correspondence is easily verified.
$\hfill$ $\qed$
\\
 
\noindent{\bf Remarks (4)}
The conclusions in the above result may fail if $\V$ is not simple, see Example 1.

\section{Irreducibility of spectral cones}
We now provide another  application of Theorem \ref{G is arcwise connected}.
Recall that a (nonempty) set $K$ is a convex cone if it is convex and $tx\in K$  for all $x\in K$ and $t\geq 0$ in $\R$. If, in addition, $K\cap -K=\{0\}$, then $K$ is said to be a pointed convex cone. We say that a  convex cone $K$ is {\it reducible} if it can be written as a sum $K=K_1+K_2$ where $K_1$ and $K_2$ are
nonzero
convex cones with $\spn(K_1)\cap \spn(K_2)=\{0\}$. A  convex cone that is not reducible is said to be {\it irreducible.}
(These concepts hold in $\Rn$ as well.) \\
We now address the irreducibility of spectral cones. In $\Rn$, spectral cones (which are just permutation invariant convex cones) can be reducible or irreducible. In fact,
$\Rnp$ ($n>1$) is reducible, while $\R_+\,{\bf 1}$ is irreducible, where ${\bf 1}$ represents the vector of ones in $\Rn$. 
 In a simple Euclidean Jordan algebra, the corresponding symmetric cone ($=$the cone of squares which is a closed convex self-dual homogeneous cone) 
is irreducible (see \cite{faraut-koranyi}, Prop. III.4.5), 
even though it  is of the form $\lambda^{-1}(\Rnp)$ with $\Rnp$  reducible for $n>1$. Below, we show that any pointed spectral cone in a simple Euclidean Jordan algebra is irreducible.
First, we recall a well-known result (\cite{hauser-guler}, slightly modified to suit our setting):
{\it If $K$ is a nonzero pointed reducible convex cone in $\V$, then there exists a unique set of  nonzero irreducible
convex cones $K_i$, $i=1,2,\ldots, r$, such that
$$K=K_1+K_2+\cdots+K_r$$
with $\spn(M)\cap \spn(N)=\{0\}$, whenever $M$ denotes the sum of  some $K_i$s and $N$  denotes the sum of the rest of the $K_i$s. Moreover, the above representation is unique up to permutation of indices.}
\\In this setting, we say that $K$ is a {\it direct sum} of $K_i$s.

\gap

\begin{theorem}
{\it Suppose $\V$ is simple. Then, every pointed spectral cone in $\V$ is irreducible.
}
\end{theorem}

\noindent{\bf Proof.}
If possible, suppose $K$ is a nonzero pointed spectral cone in $\V$ which is a direct sum of irreducible convex cones:
$$K=K_1+K_2+\cdots+K_r$$
with $r>1$.
We {\it claim that each $K_i$ is a pointed spectral cone}.
The pointedness of $K_i$ is clear, as $K$ is pointed. We show that $K_1$ is a spectral cone.  Since $\V$ is simple, because of Items $(i)$ and $(ii)$ in Proposition \ref{spectral set properties}, it is enough to show that for all $x\in K_1$,
$[x]=\{\phi(x):\phi\in G\}\subseteq K_1$.
Let $\phi\in G$. Since $K$ is a spectral cone, $\phi(K)=K$ (by Item $(iii)$ in Proposition \ref{spectral set properties}).
This implies that
$K=\phi(K_1)+\phi(K_2)+\cdots+\phi(K_r)$ is another direct sum representation in terms of irreducible convex cones. 
By the uniqueness of factors in the decomposition, $\phi(K_1)$ must be equal to some $K_j$. So,
$$\phi(K_1)\subseteq \bigcup_{j=1}^{r}\,K_j,\,\,\forall\,\,\phi\in G.$$
Now, let $0\neq x\in K_1$. Then, all the elements in $[x]$ are nonzero and
$$[x]=\{\phi(x):x\in G\}\subseteq A\cup B,$$
where $A:=K_1\,\backslash\, \{0\}$ and $B:=\cup_{2}^{r}\, K_j\,\backslash\,\{0\}.$
Since $K_1\cap \overline{B}\subseteq \spn(K_1)\cap \spn(B)=\{0\}$, we see that
$A\cap \overline{B}=\emptyset$. (Here, `overline' denotes the closure.) Similarly, $B\cap \overline{A}=\emptyset.$
So the sets $A$ and $B$ are separated \cite{rudin}. Now, by Theorem \ref{connectedness result} (or by Theorem \ref{G is arcwise connected}), $[x]$ is connected. As $[x]\subseteq A\cup B$ and $0\neq x\in A$, we must have $[x]\subseteq A$. We conclude that $[x]\subseteq K_1$. This inclusion also holds for $x=0$ so
$$[x]\subseteq K_1,\,\,\forall \,\,x\in K_1.$$
Thus, $K_1$ is a spectral cone by Item $(ii)$ in Proposition \ref{spectral set properties}.
A similar argument works for all other cones $K_i$. Now,
all the spectral cones $K$, $K_1,K_2,\ldots, K_r$ are pointed spectral cones. Hence, by Theorem 7.3 in \cite{jeong-gowda-spectral cone}, $e$ (the unit element of $\V$) or $-e$ belongs to all of them. Since $K$ is pointed, either $e$ or$-e$, but not both, can belong to all. Moreover, since the cones $K_i$ have zero as the only common element,
we conclude that $r=1$. This violates our assumption that $r>1$.
Hence, $K$ is irreducible.
$\hfill$ $\qed$

\gap

The above result gives an alternate proof of the fact that in any simple algebra, the symmetric cone is irreducible. Also, every pointed convex cone $K$ in ${\cal S}^n$ satisfying the condition
$$X\in K\Rightarrow UXU^T\in K \,\,\mbox{for all}\,\,\mbox{orthogonal matrices}\,\,U\in \R^{n\times n}$$
is irreducible.  We provide two more examples. 

\gap

\noindent{\bf Example (3)}
Let 
$n\geq 3$ and $m\in \{1,2,\ldots, n-1\}$. For each $q\in \Rn$, let $s_m(q)$ be the sum of the smallest $m$ entries of $q$.
Then, the `rearrangement cone' 
$$Q^n_m=\{q\in \Rn: s_m(q)\geq 0\},$$
is a permutation invariant pointed closed convex cone in 
$\Rn$ \cite{jeong-thesis}. Hence, by Item $(c)$ in Proposition \ref{basic properties},
$$\lambda^{-1}(Q^n_m)=\{x\in \V: \lambda_n(x)+\lambda_{n-1}(x)+\cdots+\lambda_{n-m+1}(x)\geq 0\}$$ 
is a pointed convex cone. (Note that when $m=1$, this cone is the symmetric cone of $\V$.)
By the above theorem, this cone, in a simple algebra,  
is irreducible. 

\gap

\noindent{\bf Example (4)} Let 
$n\geq 3$ and consider the following permutation invariant cone (see Example 2 in \cite{jeong-gowda-spectral cone})
$$Q:=\Big \{q\in \Rn:\,tr(q)\geq \sqrt{\frac{n}{2}}\,||q||_2\Big \}$$ where $tr(q)$ denotes the sum of all entries of $q$ and $||q||_2$ denotes the $2$-norm of $q$.
This cone is a proper cone (that is, it is a closed convex pointed cone with nonempty interior). By the above theorem, when $\V$ is a simple algebra of rank $n$, the proper spectral cone 
$$K:=\lambda^{-1}(Q)=\Big \{x\in \V:\,tr(x)\geq \sqrt{\frac{n}{2}}\,||x||_2\Big \}$$ is irreducible, where $tr(x)$ denotes the sum of all eigenvalues of $x$ and $||x||_2:=||\lambda(x)||_2$.
Using the linearity of the trace and the strict convexity of $||\cdot||_2$, it is easy to show that every boundary vector in $K$ is an extreme vector. As $n$ (the rank of $\V$) is at least $3$, such a property is false for the symmetric cone of $\V$
(take a Jordan frame $\{e_1,e_2,\ldots, e_n\}$ and consider $e_1+e_2$ which is on the boundary of the symmetric cone of $\V$, but not an extreme vector);
so, $K$ and the symmetric cone of $\V$ are not isomorphic.
This shows that in every simple algebra of rank $n\geq 3$, there is a proper (irreducible)  spectral cone that is not isomorphic to the corresponding symmetric cone.  


\begin{thebibliography}{}

\bibitem{baes} M. Baes, \emph{Convexity and differentiability properties of spectral functions in Euclidean Jordan algebras,}
Linear Algebra  Appl., 422 (2007) 664-700.

\bibitem{borwein-lewis} J. Borwein and A.S. Lewis, \emph{Convex Analysis and Nonlinear Optimization}, Springer-Verlag, New York, 2006.

\bibitem{chandrasekaran et al}
V. Chandrasekaran, P.A. Parrilo, and A.S. Willsky, \emph{Convex graph invariants}, SIAM Review, 54 (2012) 513-541.


\bibitem{daniilidis et al} 
A. Daniilidis, A. Lewis, J. Malick, and H. Sendov, \emph{Prox-regularity of spectral functions and spectral sets}, J. Convex Anal., 15 (2008) 547-560.


\bibitem{daniilidis et al 2}
A. Daniilidis, D. Drusvyatskiy, and A.S. Lewis, \emph{Orthogonal invariance and identifiability}, SIAM J. Matrix. Anal.,
35 (2014) 580-598.

\bibitem{duistermaat-kolk} J.J. Duistermaat and J.A.C. Kolk, \emph{Lie Groups}, Springer-Verlag, Berlin, 2000.

\bibitem{fan} K. Fan, \emph{On a theorem of Weyl concerning eigenvalues of linear transformations I,} Proc. Nat. Acad. Sci. U.S.A., 35 (1949) 652-655.

\bibitem{faraut-koranyi}
J. Faraut and A. Kor\'{a}nyi, \emph{Analysis on Symmetric Cones}, Clarendon Press, Oxford, 1994.

\bibitem{gowda-jeong} M.S. Gowda and J. Jeong, \emph{Commutation principles in Euclidean Jordan algebras and normal decomposition systems}, SIAM J. Optim., 27 (2017) 1390-1402.

\bibitem{hall} B.C. Hall, \emph{Lie Groups, Lie Algebras, and Representations},
Springer-Verlag Graduate Texts in Mathematics, New York, 2003.

\bibitem{hauser-guler} R. Hauser and O. G\"{u}ler, \emph{Self-scaled barrier functions on symmetric cones and their classification,}, Foundations of Comput.  Math., 2 (2002) 121-143.

\bibitem{henrion-seeger}
R. Henrion and A. Seeger, \emph{Inradius and circumradius of various convex cones arising in applications}, Set Valued Analysis, 18 (2010) 483-511.


\bibitem{iusem-seeger} A. Iusem and A. Seeger,
\emph{Angular analysis of two classes of non-polyhedral convex cones: the point of view of optimization theory}, Computational \& Appl. Math., 26 (2007) 191-214.

\bibitem{jeong-thesis} J. Jeong, \emph{Spectral sets and functions on Euclidean Jordan algebras,} PhD Thesis, University of Maryland, Baltimore County, 2017.

\bibitem{jeong-gowda-spectral cone}
J. Jeong and M.S. Gowda, \emph{Spectral cones in Euclidean Jordan algebras},
Linear Algebra Appl., 509 (2016) 286-305.

\bibitem{jeong-gowda-spectral set}
J. Jeong and M.S. Gowda, \emph{Spectral sets and functions in Euclidean Jordan algebras}, Linear Algebra Appl., 518 (2017) 31-56.

\bibitem{lewis}
A. S. Lewis, \emph{Group invariance and convex matrix analysis}, SIAM J.
 Matrix Anal., 17 (1996) 927-949.

\bibitem{lewis2}
A. S. Lewis, \emph{Convex analysis on the Hermitian matrices}, SIAM J. Optim.,
6 (1996) 164-177.

\bibitem{lewis-sendov}
A.S. Lewis and H.S. Sendov, \emph{Twice differentiable spectral functions}, SIAM J. Matrix Anal., 23 (2001) 368-386.


\bibitem{ramirez et al} H. Ramirez, A. Seeger, and D. Sossa, \emph{Commutation principle for variational problems on Euclidean Jordan algebras}, SIAM J. Optim., 23 (2013) 687-694.

\bibitem{rudin} W. Rudin, \emph{Principles of Mathematical Analysis}, McGraw-Hill, 2006. 

\bibitem{seeger}
A. Seeger, \emph{Convex analysis of spectrally defined matrix functions}, SIAM J. on Optim., 7 (1997) 679-696.

\bibitem{sossa}
D. Sossa, \emph{Euclidean Jordan algebras and variational problems under conic constraints}, PhD Thesis, University of Chile, 2014.

\bibitem{sun-sun}
D. Sun and J. Sun, \emph{L{\"o}wner's operator and spectral functions in Euclidean Jordan algebras}, Math.  Operations Res., 33 (2008) 421-445.

\bibitem{tam} T.-Y. Tam, \emph{An extension of a result of Lewis}, Electronic J. Linear Algebra, 5 (1999) 1-10. 
\end{thebibliography}
\end{document}